\documentclass{article}
\usepackage{amssymb,amsmath,amsthm}
\newtheorem{thm}{Theorem}[section]
\begin{document}
\title{About one characteristic of set of linear
Hamiltonian systems}
\author{I.~Kopshaev
\thanks{Institute of mathematics of NAS of KAZAKHSTAN, 125 Pushkina str.,
050010 Almaty, KAZAKHSTAN email: {\em{kopshaev@math.kz}}}}

\maketitle
\date
\begin{abstract}
The families of morphisms of vector fibre bundle (\cite{Mill1})
defined by the linear Hamiltonian systems of differential
equations is considered. Authors proved that the specified
families of morphisms is not saturated (\cite{Mill2}).
\end{abstract}

We consider the vector fibre bundle $(E,p,B)$ with $R^{n}$ as a
fibre and $B$ as a base (where $B$ is full metric space). On
$(E,p,B)$ has fixed some Riemannian metric (\cite{Husemoller}, P.
58--59).

Investigate the families of morphisms $\mathfrak{S}$ of linear
enlargement of dynamic system (\cite{Rahim}):
$$
(X(m), \chi(m)): (E,p,B) \to (E,p,B),
$$
($m \in N$), of the vector fibre bundle $(E,p,B)$, where
\begin{equation}
\label{kop_eq1} B=M_n, \quad E = B \times R^{2n}, \quad p = pr_1,
\end{equation}
$$ X^t(A,x) = (\chi^tA, \mathfrak{X}(t,0,A) \cdot x),
$$
$$ \chi^t A(\cdot) = A(t+(\cdot)),$$
here $M_n$ -- the space of linear Hamiltonian systems of
differential equations (\cite{Demid}), $A \in B$,
$\mathfrak{X}(\Theta, \tau, A)$ -- Cauchy matrix of the system
$$\dot{x} = J \cdot A(t) \cdot x, \quad x \in R^{2n},$$
where $J= \left(
\begin{array}{cc}
  0 & E_n \\
  -E_n & 0 \\
\end{array}%
\right)$ -- sympletic unit, $A(t)$ -- $2n \times 2n$-dimensional
symmetric matrix with $\sup\limits_{t \in R}|A(t)| < + \infty$.

\begin{thm}
The families of morphisms $\mathfrak{S}$ of the vector fibre
bundle (\ref{kop_eq1}) is not saturated.
\end{thm}

\bibliographystyle{plain}
\bibliography{nonexistence}

\end{document}